\documentclass[pdflatex,sn-basic]{sn-jnl}
\usepackage{bbold}
\usepackage{xcolor}
\usepackage{mathtools}
\jyear{2021}%

\theoremstyle{thmstyleone}

\theoremstyle{thmstyletwo}

\theoremstyle{thmstylethree}

\begin{document}

\title[Optimization via Rejection-Free Partial Neighbor Search]{Optimization via Rejection-Free Partial Neighbor Search}

\author[1]{\fnm{Sigeng} \sur{Chen}}\email{sigeng.chen@mail.utoronto.ca}

\author[1]{\fnm{Jeffrey} \spfx{S.} \sur{Rosenthal}}\email{jeff@math.toronto.edu}

\author[2]{\fnm{Aki} \sur{Dote}}\email{dote.aki@fujitsu.com}

\author[3]{\fnm{Hirotaka} \sur{Tamura}}\email{tamura.hirotaka@dxrlab.com}

\author[4]{\fnm{Ali} \sur{Sheikholeslami}}\email{ali@ece.utoronto.ca}

\affil[1]{\orgdiv{Department of Statistical Sciences}, \orgname{University of Toronto}, \orgaddress{\street{700 University Avenue}, \city{Toronto}, \postcode{M5G 1Z5}, \state{Ontario}, \country{Canada}}}

\affil[2]{\orgname{Fujitsu Ltd.}, \orgaddress{\street{4-1-1 Kamikodanaka Nakahara-ku}, \city{Kawasaki}, \postcode{211-8588}, \state{Kanagawa}, \country{Japan}}}

\affil[3]{\orgname{DXR Laboratory Inc.}, \orgaddress{\street{ 4-38-10 Takata-Nishi, Kohoku-ku}, \city{Yokohama}, \postcode{223-0066}, \state{Kanagawa}, \country{Japan}}}

\affil[4]{\orgdiv{Department of Electrical and Computer Engineering}, \orgname{University of Toronto}, \orgaddress{\street{10 King's College Road}, \city{Toronto}, \postcode{M5S 3G4}, \state{Ontario}, \country{Canada}}}

\date{\today}

\abstract{Simulated Annealing using Metropolis steps at decreasing temperatures is widely used to solve complex combinatorial optimization problems \citep{kirkpatrick1983optimization}. In order to improve its efficiency, we can use the Rejection-Free version of the Metropolis algorithm, which avoids the inefficiency of rejections by considering all the neighbors at every step \citep{rosenthal2021jump}. As a solution to avoid the algorithm from becoming stuck in local extreme areas, we propose an enhanced version of Rejection-Free called Partial Neighbor Search (PNS), which only considers random parts of the neighbors while applying Rejection-Free. We demonstrate the superior performance of the Rejection-Free PNS algorithm by applying these methods to several examples, such as the QUBO question, the Knapsack problem, the 3R3XOR problem, and the quadratic programming.}


\maketitle
\section{Introduction} \label{sec-introduction}

Optimization is the cornerstone of many areas, and it plays a crucial role in finding feasible solutions to real-life problems, from mathematical programming to operations research, economics, management science, business, medicine, life science, and artificial intelligence \citep{floudas2008encyclopedia}. Prior to the invention of linear and integer programming in the 1950s, optimization was characterized by several independent topics, such as optimum assignment, the shortest spanning tree, transportation, and the traveling salesman problem, which were then united into one framework \citep{schrijver2005history}. Today, combinatorial optimization plays an important role in research because most of its problems originate from practice and are dealt with on a daily basis \citep{schrijver2005history}. The process of finding a feasible solution to some complex combinatorial optimization problems may take a considerable amount of time. In particular, no algorithm for NP-hard problems can guarantee that the optimal state of the problem will be found within a limitation governed by a polynomial based on the input length \citep{garey1974some}.

In general, metaheuristics are algorithmic frameworks that are often nature-inspired and are used to solve complex optimization problems \citep{bianchi2009survey} by arriving at a feasible solution, regardless of whether it is optimal. The Simulated Annealing algorithm \citep{kirkpatrick1983optimization}, based on the Metropolis algorithm \citep{Metropolis1953} at decreasing temperatures, is a typical method of this kind. The Simulated Annealing algorithm, however, may be inefficient with respect to rejections. In order to improve the performance of Simulated Annealing, we adopt the Rejection-Free algorithm for sampling \citep{rosenthal2021jump} into an optimization version. Additionally, Rejection-Free may experience inefficiency when it enters local extreme areas. Therefore, we propose another algorithm based on the Rejection-Free algorithm called Partial Neighbor Search (PNS) in order to further enhance its efficiency.

Even when applied to a single-core implementation, Rejection-Free and PNS are more efficient in many optimization problems than Simulated Annealing. Moreover, the implementation of these algorithms can also be carried out through parallelism in order to increase efficiency even further. It is possible to use processors designed for general purposes, such as Intel and AMD cores, for parallel computing to accelerate the algorithm to some extent. However, these chips were not built for parallel computing, and off-chip communication significantly slows the data transfer rate to and from the cores \citep{sodan2010parallelism}. On the other hand, parallelism hardware designed specifically for MCMC trials has been proposed. For example, the second generation of Fujitsu Digital Annealer uses a dedicated processor called Digital Annealing Unit (DAU) \citep{9045100} to achieve high speed. This dedicated processor is designed to minimize communication overhead in arithmetic circuitry and with memory. It is possible to achieve 100x to 10,000x speedups by combining Rejection-Free and PNS with such parallelism hardware \citep{sheikholeslami2021power}.

We next review the Simulated Annealing algorithm, the Metropolis algorithm, and the Rejection-Free algorithm for sampling. Following that, Section~\ref{sec-rejectionfree} describes how to use the Rejection-Free algorithm to solve optimization problems. Our next point is that the local maximum may lead to another kind of inefficiency for Rejection-Free, and Section~\ref{sec-pns} introduces our Partial Neighbor Search (PNS) algorithm for optimization, which considers just subsets of neighbor states for possible moves. In Section~\ref{sec-qubo}, we demonstrate how PNS can be applied to quadratic unconstrained binary optimization (QUBO) questions and its effectiveness in solving them. We then discuss why this improvement occurs (Section~\ref{sec-understand}), and how its subsets of partial neighbors should be chosen (Section~\ref{sec-optimal}), as well as its relation to the Tabu Search algorithm (Section~\ref{sec-tabu}). Moreover, we present several other examples, such as the Knapsack problem (Section~\ref{sec-knapsack}) and the 3R3XOR problem (Section~\ref{sec-xor}), to illustrate the advantages of the PNS algorithm in discrete optimization problems. Furthermore, Section~\ref{sec-continuous} illustrates another advantage of PNS over Rejection-Free by providing a continuous optimization example known as quadratic programming. PNS can easily be adapted to the general state space by selecting only a finite subset, and it outperforms Simulated Annealing, whereas Rejection-Free cannot be applied in this case due to the need to consider all neighbors at each step.

\subsection{Background on Simulated Annealing for optimization} \label{subsec-sa}

Simulated Annealing, as introduced by \cite{kirkpatrick1983optimization}, is widely used to solve combinatorial optimization problems, such as approximating the optimal values of functions with many variables \citep{rutenbar1989simulated}. Although there is no guarantee that this algorithm will provide an optimal solution, it is capable of providing reasonable, feasible solutions quickly \citep{albright2007introduction}. Simulated Annealing contains the following essential elements \citep{bertsimas1993simulated}:

\begin{enumerate}
    \item A state space $\mathcal{S}$.
    \item A real-valued target distribution $\pi$ on $\mathcal{S}$. The ultimate goal for the Simulated Annealing is to find $Y \in \mathcal{S}$ such that $\pi(Y) > \pi(X)$, $\forall X \in \mathcal{S}$. However, for many circumstances, a good feasible solution is acceptable.
    \item $\forall X \in S$, $\exists$ a proposal distribution $\mathcal{Q}(X, \cdot)$ where $\int_{Y \in \mathcal{S} \backslash \{X\}} \mathcal{Q}(X, Y) = 1$.
    \item $\forall X \in \mathcal{S}$, $\exists$ $\mathcal{N} (X) = \{Y \in \mathcal{S} \mid \mathcal{Q}(X, Y) > 0\} \subset \mathcal{S} \backslash \{X\}$, called the neighbors of $X$.
    \item A non-increasing function $T: \mathbb{N} \to (0, \infty)$, called the Cooling Schedule. $T(k)$ is called the temperature at step $k \in \mathbb{N}$.
    \item An initial State $X_0 \in \mathcal{S}$.
\end{enumerate}

With the above elements, the Simulated Annealing algorithm, which consists of a discrete time-inhomogeneous Markov Chain $\{X_k\}_{k = 0}^K$ can be generated by Algorithm~\ref{alg-sa}. Algorithm~\ref{alg-sa} is designed to converge to states $X_k$ with nearly-maximal values of $\pi(X_k)$, though that is not guaranteed. Note that the algorithm can also be formulated using log values for better numerical stability.

\begin{algorithm}
\caption{Simulated Annealing}\label{alg-sa}
\begin{algorithmic}
\State initialize $X_0$
\For{$k$ in $1$ to $K$}
    \State random $Y \in \mathcal{N}(X_{k-1})$ based on $\mathcal{Q}(X_{k-1}, \cdot)$
    \State random $U_k \sim \text{Uniform}(0, 1)$
    \If{$U_k < \Large[\frac{\pi(Y)}{\pi(X_{k-1})}\Large]^{{1 / T(k)}}$} 
    \State \Comment{accept with probability $\min \Big{\{} 1, \big{[}\frac{\pi(Y_k)}{\pi(X_{k-1})}\big{]}^{{1 / T(k)}} \Big{\}} $}
    \State $X_{k} = Y$ \Comment{accept and move to state $Y$}
    \Else \State $X_{k} = X_{k-1}$ \Comment{reject and stay at $X_{k-1}$}
    \EndIf 
\EndFor
\end{algorithmic}
\end{algorithm}

\subsection{Background on Metropolis-Hastings algorithm} \label{subsec-metropolis}

The above Simulated Annealing algorithm is designed based on the Metropolis algorithm \citep{Metropolis1953}. Among all the Monte Carlo algorithms, the Metropolis algorithm has been the most successful and influential \citep{beichl2000metropolis}. It is designed to generate a Markov chain that converges to a given target distribution $\pi$ on a state space $\mathcal{S}$. As a generalization of the Metropolis algorithm, the Metropolis-Hastings(M-H) algorithm includes the possibility of a non-symmetric proposal distribution $\mathcal{Q}$ \citep{hitchcock2003history}. The M-H algorithm is described in Algorithm~\ref{alg-metropolis}.
 
\begin{algorithm}
\caption{the Metropolis-Hastings algorithm}\label{alg-metropolis}
\begin{algorithmic}
\State initialize $X_0$
\For{$k$ in $1$ to $K$}
    \State random $Y \in \mathcal{N}(X_{k-1})$ based on $\mathcal{Q}(X_{k-1}, \cdot)$
    \State random $U_k \sim \text{Uniform}(0, 1)$
    \If{$U_k < \frac{\pi(Y) \mathcal{Q}(Y, X_{k-1})}{\pi(X_{k-1}) \mathcal{Q}(X_{k-1}, Y)}$}     
    \State \Comment{accept with probability $\min \Big{\{} 1,  \frac{\pi(Y) \mathcal{Q}(Y, X_{k-1})}{\pi(X_{k-1}) \mathcal{Q}(X_{k-1}, Y)} \Big{\}} $}
    \State $X_{k} \gets Y$ \Comment{accept and move to state $Y$}
    \Else \State $X_{k} \gets X_{k-1}$ \Comment{reject and stay at $X_{k-1}$}
    \EndIf 
\EndFor
\end{algorithmic}
\end{algorithm}

Algorithm~\ref{alg-metropolis} ensures the Markov chain $\{X_0, X_1, X_2, \dots, X_K\}$ has $\pi$ as stationary distribution. It follows (assuming irreducibility) that the expected value $E_\pi(h)$ of a functional $h: \mathcal{S} \to \mathbb{R}$ with respect to $\pi$ can be estimated by $\frac{1}{M} \sum_{i=1}^M h(X_i)$ for sufficiently large run length $M$. Although the M-H algorithm and Simulated Annealing are designed for different purposes, regarding the implementation, the Cooling Schedule is the only difference between them. Thus, both Simulated Annealing and the M-H algorithm may face inefficiencies from the rejections \citep{rosenthal2021jump}.

\subsection{Background on Rejection-Free algorithm for sampling} \label{subsec-rfsampling}

Rejections in both Simulated Annealing and the M-H algorithm could be a problem. In Algorithm~\ref{alg-sa}, if $U_k \ge \Large[\frac{\pi(Y)}{\pi(X_k)}\Large]^{{1 / T(k)}}$, then we will remain at the current state, even though we have spent time in proposing a state, computing a ratio of target probabilities, generating a random variable $U_k$, and deciding not to accept the proposal. Such inefficiencies could happen frequently and are considered a necessary evil of Simulated Annealing and the M-H algorithm. However, we can compute all potential acceptance probabilities at once to allow for the possibility of skipping these rejection steps \citep{rosenthal2021jump}. By taking out the inefficiencies of rejections in both algorithms, the Rejection-Free algorithm can lead to significant speedup. 

Before introducing Rejection-Free, we need to introduce the jump chain first. Given a run $\{X_k\}$ of a Markov chain, we define the jump chain to be $\{J_k, M_k\}$, where $\{J_k\}$ represents the same chain as $\{X_k\}$ except omitting any immediately repeated states, and the Multiplicity List $\{M_k\}$ is used to count the number of times the original chain remains at the same state.

For example, if the original chain is 
\begin{equation}
    \{X_k\} = \{a, b, b, b, a, a, c, c, c, c, d, d, a, \dots\},
\end{equation}
then the jump chain would be 
\begin{equation}
\{J_k\}= \{a, b, a, c, d, a, \dots\},
\end{equation}
with the corresponding multiplicity list being
\begin{equation}
\{M_k\} = \{1, 3, 2, 4, 2, 1, \dots\}.
\end{equation}
The jump chain $\{J_k, M_k\}$ itself is a Markov chain, with transition probabilities $\hat{P}(y \mid x)$ specified by 
\begin{equation}
\begin{aligned}
    \hat{P}(x \mid x) & \coloneqq 0\\
    \forall y \ne x \mbox{,    } \hat{P}(y \mid x) & \coloneqq P[J_{k+1} = y \mid J_k = x] = \frac{P(y \mid x)}{\sum_{z \ne x}P(z \mid x)} 
\end{aligned}
\end{equation}
Moreover, the conditional distribution of $\{M_k\}$ given $\{J_k\}$ is equal to the distribution of $1 + G$ where $G$ is a geometric random variable with success probability $p = 1 - P(x \mid x) = \sum_{z \ne x}P(z \mid x)$; see \cite{rosenthal2021jump}. 

Given the above properties for the Jump chain, the Rejection-Free algorithm can be used for sampling as described by Algorithm~\ref{alg-rfsampling}. Algorithm \ref{alg-rfsampling} only works for the discrete cases where all states have at most finite neighbors. Theorem 13 in \cite{rosenthal2021jump} extended the Rejection-Free to general state space, and we will discuss more by a continuous optimization question in Section~\ref{sec-continuous}.

\begin{algorithm}
\caption{Rejection-Free for Sampling (Discrete Cases)}\label{alg-rfsampling}
\begin{algorithmic}
\State initialize $J_0$
\For{$k$ in $1$ to $K$}
    \State $p \gets 0$  \Comment{p is used to record the success probability for $M_{k-1}$}
    \For{Y in $\mathcal{N}(J_{k-1})$}  \Comment{only works for finite neighbors}
        \State calculate $q(Y) = \mathcal{Q}(Y, J_{k-1})\min\{1, \frac{\pi(y)}{\pi(J_{k-1})}\}$ 
        \State \Comment{the transition prob. from $J_{k-1}$ to $Y$}
        \State $p \gets p + q(Y)$ \Comment{$p = \sum_{z \ne x}P(z \mid x)$}
    \EndFor
    \State choose $J_{k} \in \mathcal{N}(J_{k-1})$ such that $\hat{P}(J_{k} = Y \mid J_{k-1}) \propto q(Y)$
    \State \Comment{choose the next jump chain state}
    \State calculate $M_{k-1} = 1 + G$ where $G \sim \text{Geom}(p)$
    \State \Comment{multiplicity list for current state}
\EndFor
\end{algorithmic}
\end{algorithm}

Algorithm~\ref{alg-rfsampling} ensures (assuming irreducibility) that the expected value $E_\pi(h)$ of a functional $h: \mathcal{S} \to \mathbb{R}$ with respect to $\pi$ can be estimated by $\frac{\sum_{k=1}^K M_k \, h(J_k)}{\sum_{k=1}^K M_k}$ for sufficiently large run length $K$, while avoiding any rejections. Rejection-Free can lead to great speedup in examples where the Metropolis algorithm frequently rejects \citep{rosenthal2021jump}.

\section{Rejection-Free algorithm for optimization} \label{sec-rejectionfree}

In addition to sampling, the above Rejection-Free algorithm can also be applied to optimization problems. Given a set $\mathcal{S}$ and a real-valued target distribution $\pi$ on the set $\mathcal{S}$, we can use the Rejection-Free algorithm to find $X \in \mathcal{S}$ that maximizes $\pi(X)$ by Algorithm~\ref{alg-rf}. Algorithm~\ref{alg-rf} is again designed to converge to states $X_k$ with nearly-maximal values of $\pi(X_k)$, with greater efficiency by avoiding rejections, though that is again not guaranteed. Although the purpose of sampling and optimization are different, regarding the implementation, Rejection-Free for optimization is only different from Rejection-Free for sampling by getting rid of the multiplicity list $\{M_k\}$.

\begin{algorithm}
\caption{Rejection-Free for Optimization (Discrete Cases)}\label{alg-rf}
\begin{algorithmic}
\State initialize $J_0$
\For{$k$ in $1$ to $K$}
    \For{$Y \in \mathcal{N}(J_{k-1})$}
    \State \Comment{only works for finite neighbors}
        \State calculate $q(Y) = \mathcal{Q}(Y, J_{k-1}) \min\{1, [\frac{\pi(Y)}{\pi(J_k)}]^{\frac{1}{T(k)}}\}$  
        \State \Comment{the transition prob. from $J_{k-1}$ to $Y$}
    \EndFor
    \State choose $J_{k} \in \mathcal{N}(J_{k-1})$ such that $\hat{P}(J_{k} = Y \mid J_{k-1}) \propto q(Y)$ \State \Comment{choose the next jump chain State}
\EndFor
\end{algorithmic}
\end{algorithm}

\begin{figure}
    \centering
    \includegraphics[width=\textwidth]{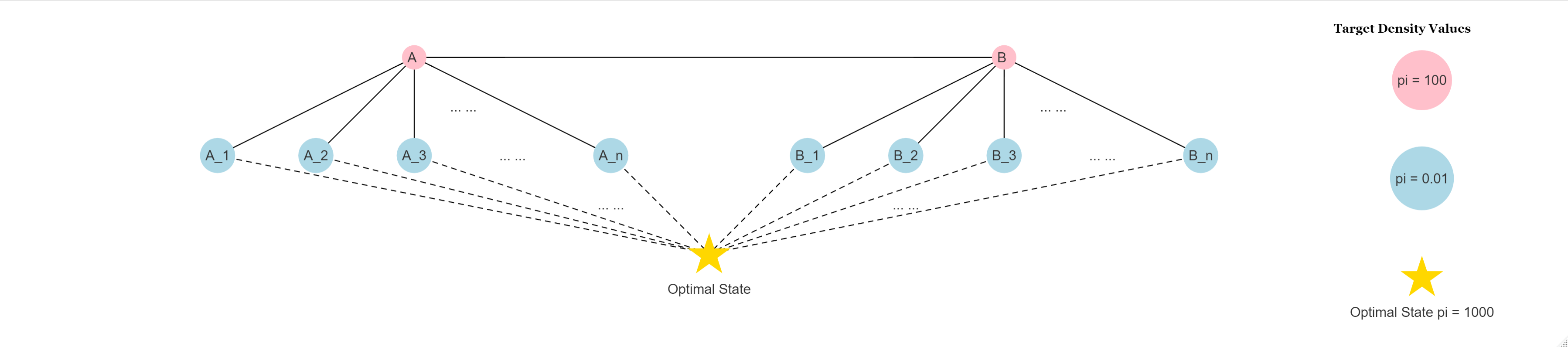}
    \caption{Illustration of the local maximum area in an optimization problem where both Simulated Annealing and Rejection-Free may get stuck. The target distribution $\pi$ has the following function values: $\pi(A) = \pi(B) = 100$, $\pi(A_1) = \pi(A_2) = \dots = \pi(A_n) = \pi(B_1) = \pi(B_2) = \dots = \pi(B_n) = 0.01$.}
    \label{fig-stuckexample}
\end{figure}

Although the Rejection-Free algorithm for optimization can help reduce the inefficiency of rejections, local maximum areas of $\pi$ can still be a problem. For example, we want to find $X \in S$, which maximizes $\pi(X)$ from a state space starting at state $A$ in Figure~\ref{fig-stuckexample}. Here, we use a uniform proposal distribution $\mathcal{Q}$ on the neighbor sets $\mathcal{N}$ as shown in Figure \ref{fig-stuckexample}. We will have many rejections if we constantly use Simulated Annealing with $T \equiv 1$. Note that, $\pi(A_1) = \pi(A_2) = \dots = \pi(A_n) = 0.01$ while $\pi(A) = \pi(B) = 100$. The probability of escaping from A is $\frac{1}{n+1} + \frac{n}{n+1} \times \frac{1}{10000}$, where $\frac{1}{n+1}$ represents the probability of moving from state A to state B, and $\frac{n}{n+1} \times \frac{1}{10000}$ is the probability of moving from state A to $A_1, A_2, \dots, A_n$. Cooling Schedules can help reduce the probability of rejection at the beginning of Simulated Annealing since $T$ should be large at the beginning. However, as we move on in Simulated Annealing, we will be more and more likely to be trapped by local maximum areas like this. The Rejection-Free algorithm for optimization can produce some speedup in this case, but the Rejection-Free chain will still be stuck by the local maximum area $\{\pi(A), \pi(B)\}$. If $n$, the number of other neighbors for $A$ and $B$, is small, this chain will be switching between A and B for a really long time, since 
\begin{equation}
\begin{aligned}
    \hat{P}(J_{1} & = B \mid J_0 = A) = \frac{\min\{1, \frac{\pi(B)}{\pi(A)}\}}{\sum_{z \ne A}\min\{1, \frac{\pi(z)}{\pi(A)}\}} = \frac{1}{1 + 0.0001 \times n} \approx 1 \\
    \hat{P}(J_{1} & = A \mid J_0 = B) \approx 1.
\end{aligned}
\end{equation}
To help our Markov chain escape from those local maximums in optimization, we propose another method called Partial Neighbor Search based on the Rejection-Free algorithm. 

\section{Proposed Search Algorithm: Partial Neighbor Search} \label{sec-pns}

Partial Neighbor Search (PNS) is an algorithm based on the Rejection-Free, also designed as a Markov chain used for optimization as described in Algorithm~\ref{alg-pns}. Algorithm~\ref{alg-pns} is again designed to converge to states $X_k$ with nearly-maximal values of $\pi(X_k)$, with greater efficiency by avoiding both rejections and traps in local maximum areas.

\begin{algorithm}
\caption{Partial Neighbor Search}\label{alg-pns}
\begin{algorithmic}
\State initialize $J_0$
\For{$k$ in $1$ to $K$}
    \State pick $\mathcal{N}_k(J_{k-1}) \subset \mathcal{N}(J_{k-1})$ $(\star)$ 
    \For{$Y \in \mathcal{N}_k(J_{k-1})$} \Comment{Only neighbors in $\mathcal{N}_k$ will be considered}
        \State calculate $q(Y) =  \mathcal{Q}(Y, J_{k-1}) \min\{1, [\frac{\pi(Y)}{\pi(J_k)}]^{\frac{1}{T(k)}}\}$
        \State \Comment{the transition prob. from $J_{k-1}$ to $Y$}
    \EndFor
    \State choose $J_{k} \in \mathcal{N}_k(J_{k-1})$ such that $\hat{P}(J_{k} = Y \mid J_{k-1}) \propto q(Y)$ 
    \State \Comment{choose the next Jump Chain State}
\EndFor
\end{algorithmic}
\end{algorithm}

The $(\star)$ step in Algorithm~\ref{alg-pns} is the key of PNS. At this step, $\mathcal{N}_k(J_{k-1})$ could be random $50\%$ of the elements from $\mathcal{N}(J_{k-1})$. In Section~\ref{sec-optimal}, we will explore many other choices for the $(\star)$ step to figuring out the best strategy. Moreover, for continuous cases, PNS can be applied, and we only need to ensure the Partial Neighbor Sets $\mathcal {N}_k$ are always finite, $\forall k$. On the other hand, Algorithms~\ref{alg-rfsampling} and Algorithm~\ref{alg-rf} for Rejection-Free only work for discrete cases where the number of neighbors for all states must be finite, and we will illustrate these by an optimization example in continuous cases in Section~\ref{sec-continuous}. 

The motivation for PNS is simple: we have a better chance of escaping from the local maximum area if we force the algorithm to avoid some neighbors randomly. For example, in Figure~\ref{fig-stuckexample}, if we only consider half of the neighbors at state $A$, then we may disregard state $B$ with probability $50\%$, then we have a probability of at least $50\%$ of selecting a state from $\{A_1, A_2, \dots, A_n\}$ as our next state in the PNS chain. If this occurs, we are more likely to escape from the local maximum area $\{ \pi(A), \pi(B)\}$.

\section{Application to the QUBO question} \label{sec-qubo}

The quadratic unconstrained binary optimization (QUBO) has gained increasing attention in the field of combinatorial optimization due to its wide range of applications in finance and economics to machine learning \citep{kochenberger2014unconstrained}. The QUBO problem is known to be NP-hard \citep{glover2018tutorial}, so it is common to use Simulated Annealing to find the optimal or workable solution. This problem can now be addressed using our PNS algorithm. (Additional applications are in Sections~\ref{sec-knapsack}, Section~\ref{sec-xor}, and Section~\ref{sec-continuous} below.)

For a given $N$ by $N$ matrix $Q$ (usually upper triangular), the QUBO question aims to find 
\begin{equation}
    \arg \max x^T Q x \mbox{, where } x \in \{0, 1\}^N
\end{equation}
(Sometimes $\arg\min$ is used in place of $\arg\max$, which is equivalent to taking the negative of $Q$, so for simplicity, we focus on the $\arg\max$ version here.)

As part of our algorithm, we use a uniform proposal distribution among all neighbors where the neighbors are defined as binary vectors with Hamming distance 1. That is, $\mathcal{Q}(X, Y) = \frac{1}{N}$ for $\forall Y \in \mathcal{N}(X)$,  where $Y \in \mathcal{N}(X) \iff \lvert X - Y \rvert = \sum_{i=1}^N \lvert X_i - Y_i \rvert = 1$, $\forall X, Y \in \{0, 1\}^N$.  We randomly choose half of the neighbors at each step of PNS, which means we only consider a random subset $\mathcal{N}_K(x) \subset \mathcal{N}(x)$ whose cardinality is $\lvert \mathcal{N}_k(X) \rvert = \frac{1}{2} \lvert \mathcal{N}(X) \rvert =  \frac{1}{2}N$ for $\forall X \in \{0, 1\}^N$. In addition, the target distribution $\pi(x) = \exp\{x^T Q x\}$, since we need the target distribution to be positive all time to use the Cooling Schedule, and maximizing $x^T Q x$ is the same as maximizing $\exp\{x^T Q x\}$. Furthermore, $T(k)$ represents the temperature at step $k$ for the cooling schedule here.

We compare Simulated Annealing, Rejection-Free for Optimization, and PNS in 1000 simulation runs. We randomly generate a 200 by 200 upper triangular as the QUBO matrix $Q$. The non-zero elements from $Q$ were generated randomly by $Q_{i, j} \sim \text{Normal}(0, 100^2)\mbox{, } \forall i \le j$.

\begin{figure}
    \centering
    \includegraphics[width=\textwidth]{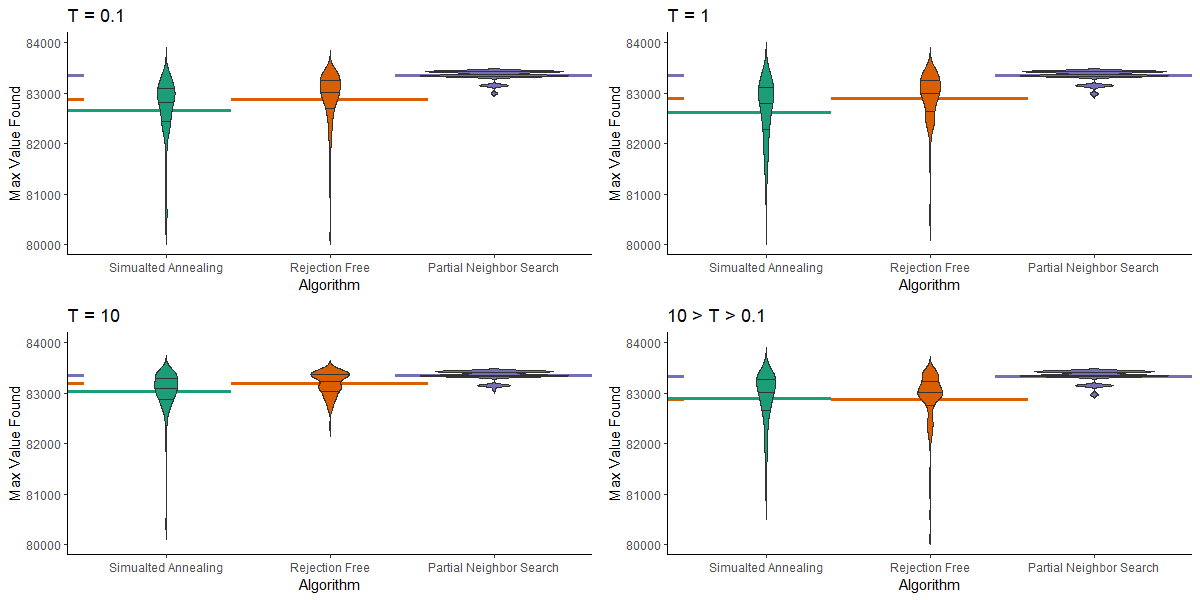} \\
    \caption{Comparison of Simulated Annealing, Rejection-Free, and PNS in terms of the highest (log) target distribution value $\log \pi(x) = x^T Q x$ being found, for a random upper triangular QUBO matrix $Q$ where the non-zero elements are generated by $Q_{i, j} \sim N(0, 100^2)$. Four different cooling schedules where $T(k) = 0.1$, $1$ and $10$ constantly, and $T(k)$ being geometric from 10 to 0.1 are used here. The number of iterations for Simulated Annealing is 200,000, and the numbers of iterations for Rejection-Free and PNS are 1000. The three black lines inside the violin plots are $25\%$, $50\%$, and $75\%$ quantile lines. The colored segments represent the mean values.}
    \label{fig_compare_sa_rf_pns}
\end{figure}

The result for the simulation is shown in Figure \ref{fig_compare_sa_rf_pns}. Here, we used a violin plot to summarize the results. The violin plot uses the information available from local density estimates and the basic summary statistics to provide a useful tool for data analysis and exploration \citep{hintze1998violin}. The violin plot combines two density traces on both sides and three quantile lines ($25\%$, $50\%$, and $75\%$) to reveal the data structure. In addition, we added a long segment of the bottom layer as the mean for the values. We also added a short segment on the y-axis to help compare the mean values. 

From Figure \ref{fig_compare_sa_rf_pns}, we can see that the PNS is always the best in all four different cooling schedules. Note that the number of iterations used for Simulated Annealing is $200,000$ for Simulated Annealing while they are $1000$ for both Rejection-Free and PNS. We used these many iterations because we need to consider $200$ neighbors at each iteration in Rejection-Free, while we only need to consider one neighbor for each iteration in Simulated Annealing. If we proceed with all three algorithms on a single-core machine, the run time of a single simulation run for simulated Annealing is about 20 seconds; the run time for Rejection-Free is about 10 seconds; the run time for PNS is only 5 seconds. In addition, parallelism in computer hardware can increase the speed of both Rejection-Free and PNS by distributing the calculation of the transition probabilities for different neighbors onto different cores \citep{rosenthal2021jump}. Besides that, we can also use multiple replicas at different temperatures, such as in parallel tempering, or deploy a population of replicas at the same temperature \citep{sheikholeslami2021power}. Combining these methods by parallelism can yield 100x to 10,000x speedups for both Rejection-Free and PNS \citep{sheikholeslami2021power}.

In the above example, the improvement in the efficiency of Rejection-Free is not hard to understand. The performance of PNS is somehow counter-intuitive. Compared to Rejection-Free, why would we get a better result by considering fewer neighbors at each step? To illustrate how PNS works, we can look closely at the Markov chains generated in the above example.

\section{Understanding the improvement of Partial Neighbor Search} \label{sec-understand}

\begin{figure} 
    \centering
    \includegraphics[width=\textwidth]{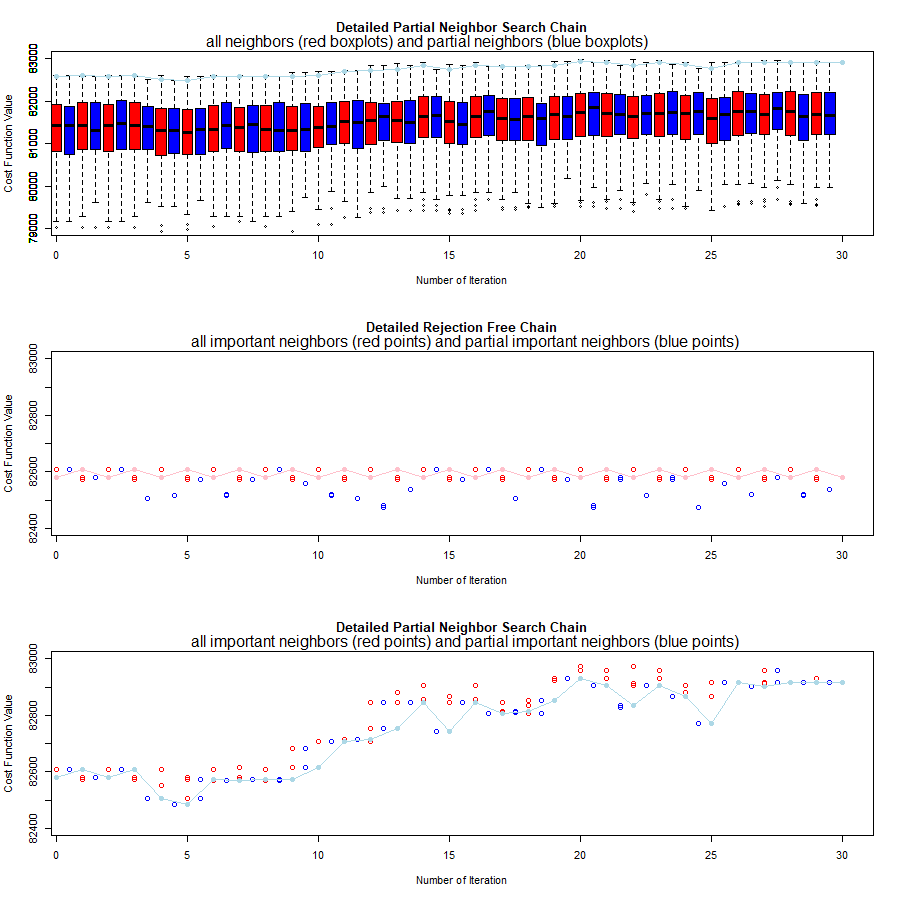}
    \caption{The detailed Markov Chains from Rejection-Free (the pink chain in the second plot) and PNS (the light blue chain in the first and the third plot). The red box plots in the first plot represent the target distribution values for all neighbors, and the blue box plots represent the partial neighbors. Most of these values are useless because they are too small to be picked by the Markov chain. The second and the third plots only show the important neighbors, defined as those whose transition probability is larger than $\exp\{-10\}$ times the highest transition probability among all neighbors. Here, red points represent all important neighbors, and blue points mean important neighbors of a random subset of all neighbors used for PNS. The Rejection-Free Chain switches between three local maximum states all the time while the PNS chain escapes from the local maximum area after five iterations.}
    \label{fig-detailedchain}
\end{figure}

In this section, we found a local maximum area for the target distribution $\pi$ purposefully in the previous QUBO example in Section~\ref{sec-qubo} by looking at the final results from the simulation runs from the previous section. Many Rejection-Free chains stopped at this local maximum area after 1000 iterations. For this local maximum area, the target distribution value is around 82600, and this local maximum area contains three states whose target distribution values are much larger than all their other neighbors. Thus, this local maximum can trap the regular Rejection-Free chain for a long time, just like the cases we mentioned in Figure~\ref{fig-stuckexample}. We can plot the Markov chains by PNS with the target distribution values for all the neighbors by Rejection-Free and the random subset of neighbors by PNS in the form of boxplots. The boxplot of the first 30 steps from the first simulation in PNS is shown in the first plot in Figure \ref{fig-detailedchain}

From the first plot in Figure \ref{fig-detailedchain}, most of the target distribution values within the boxplot are not useful since they are too small to be picked by the algorithm. Therefore, we only need to consider the important neighbors likely to be chosen. Firstly, for each state $J_k$ in the Markov Chain, we find the max value among all the transition probabilities, and we define the important neighbors to be those neighbors whose transition probability is larger than $\exp\{-10\}$ times the highest transition probability among all neighbors. That is, for each $J_k$ from the chain, we find $q(Y_0) = \max\{q(Y) \mid Y \in \mathcal{N}(J_{k})\}$, and then we define $\{Y \mid Y \in \mathcal{N}(J_k), q(Y) > \exp\{-10\} \times q(Y_0)\}$ to be important neighbors for $J_k$. This time, we only have several important neighbors at each step. Thus, we used points instead of boxplots to show the important neighbors. The result from Rejection-Free and PNS is also shown in Figure~\ref{fig-detailedchain}.

From the second plot in Figure~\ref{fig-detailedchain}, the red dots represent the important neighbors, and the pink line means the Rejection-Free chain. We can see that this local maximum area of three states can easily trap the Rejection-Free chains because their target distribution values are much higher than all other neighbors. Thus, the important neighbors for any of these three states are only the remaining two, and the Rejection-Free chain will be switching between these three for a long time. At the same time, the blue dots in the second plot represent the important neighbors if we start to do PNS from that state. Although we did not apply PNS in the second plot, we still put the random subset for PNS there as a comparison. From the blue dots in the second plot, we can say that if we perform PNS, then the Markov chain can escape from this local maximum area faster since some groups of the blue dots do not contain any of these three states with high target distribution values. 

On the other hand, the third plot in Figure~\ref{fig-detailedchain} shows that the PNS chain (blue line) escapes from this local maximum area within five steps. Again, the blue dots represent the important neighbor from PNS, and the red dots represent the important neighbor if we start to perform Rejection-Free from that step. For each step of PNS within the local maximum area of three states, the Markov Chain has the probability of $25\%$ to include neither of the remaining neighbors from the three states. Thus, PNS helped the Markov chain to escape from this local maximum area. In addition, in the middle part of the PNS chain, when the target distribution value of the PNS chain is increasing, we usually have more than one important neighbor. For example, if we have three important neighbors, we only have $12.5\%$ for considering none of them by PNS.

Thus, the PNS is better than Rejection-Free because the PNS performs much better than the Rejection-Free algorithm when the local maximum areas trap the Markov chain. On the other hand, PNS is not much worse than Rejection-Free when the Markov chain is increasing with respect to the target distribution value.

This section uses $50\%$ random partial neighbors for each step. We have many other choices, and we will consider and compare these choices in the next section.

\section{Optimal subset choice for Partial Neighbor Search} \label{sec-optimal}

We formally define the way of choosing Partial Neighbors Sets. Before we start the Markov chain, we need to define a proposal distribution $\mathcal{Q}$ and corresponding neighbor sets $\mathcal{N}(X) := \{Y \in S \mid \mathcal{Q}(X, Y) > 0\}$. Partial Neighbor Sets $\mathcal{N}_k$ means any set satisfies the following conditions:
\begin{enumerate}
    \item $\mathcal{N}_k: \mathcal{S} \to P(\mathcal{S})$, where $\mathcal{S}$ is the state space, and $P(\mathcal{S})$ is the power set of $\mathcal{S}$.
    \item $\mathcal{N}_k(X) \subset \mathcal{N}(X)$, $\forall X \in \mathcal{S}$.
    \item $Y \in \mathcal{N}_k(X)  \iff X \in \mathcal{N}_k(Y)$, $\forall X, Y \in \mathcal{S}$.
    
    \item Define $\mathcal{Q}_k(X, Y): \mathcal{S} \times \mathcal{S} \to \mathbb{R}$ be the corresponding partial proposal distribution where $\mathcal{Q}_k(X, Y) \propto \mathcal{Q}(X, Y)$ for $Y \in \mathcal{N}_k(X)$ and $\mathcal{Q}_k(X, Y) = 0$ otherwise. 
    \item Define the Partial Neighbor Weight $t \coloneqq \int_{Y \in \mathcal{N}_k(X)} \mathcal{Q}(X, y) d y$. Note that if we want to ensure the reversibility of the Markov chain, then we have to make sure the Partial Neighbor Weight is a constant.
\end{enumerate}
Usually, we want to pick $\mathcal{N}_i$ such that $\lvert \mathcal{N}_k(X) \rvert < \lvert \mathcal{N}(X) \rvert$ to perform proper PNS. In addition, to ensure irreducibility, we need to make sure $\cup_{i=0}^{K-1} \mathcal{N}_k(X)= \mathcal{N}(X)$ for all $X \in \mathcal{S}$. 

Here, we compare the four different ways to choose the proposal distribution $\{\mathcal{Q}_k, \mathcal{N}_k\}$ for PNS in the $(\star)$ step in Algorithm \ref{alg-pns}:

\begin{itemize}
    \item Method A (random subset every step): The Partial Neighbor Sets $\mathcal{N}_k$ are randomized for every step, where $\lvert \mathcal{N}_k(X) \rvert = \frac{1}{2} \times \lvert \mathcal{N}(X) \rvert$. $\mathcal{Q}_k$'s are defined accordingly.

    \item Method B (random subset every 10 steps): The Partial Neighbor Sets $\mathcal{N}_k$ are randomized for once 10 steps, where $\lvert \mathcal{N}_k(X) \rvert = \frac{1}{2} \times \lvert \mathcal{N}(X) \rvert$. That is, $\mathcal{N}_{10 \times k + 1}= \mathcal{N}_{10 \times k + 2}= \dots = \mathcal{N}_{10 \times k + 10}$ for $\forall k \in \mathbb{N}$. $\mathcal{Q}_k$'s are defined accordingly.

    \item Method C (systematic subset every step): Before we start our Markov Chain, we define two Partial Neighbor Sets $\mathcal{N}_1$ and $\mathcal{N}_2$, where $\lvert \mathcal{N}_1(X) \rvert = \lvert \mathcal{N}_2(X)\rvert = \frac{1}{2} \times \lvert \mathcal{N}(X) \rvert$, $\mathcal{N}_1(X) \cap \mathcal{N}_2(X) = \emptyset$. For step $k$ of the Markov chain, we only randomly generate $r_k \in \{1, 2\}$, and apply $\mathcal{N}_{r_k}$ for step $k$. $\mathcal{Q}_1$ and $\mathcal{Q}_2$  are defined accordingly.

    \item Method D  (systematic subset every 10 steps): Before we start our Markov Chain, we define two Partial Neighbor Sets $\mathcal{N}_1$ and $\mathcal{N}_2$, where $\lvert \mathcal{N}_1(X) \rvert = \lvert \mathcal{N}_2(X) \rvert = \frac{1}{2} \times \lvert \mathcal{N}(X)\rvert$, $\mathcal{N}_1(X) \cap \mathcal{N}_2(X) = \emptyset$. For every ten steps of the Markov chain, we only randomly generate $r_k \in \{1, 2\}$ and apply $\mathcal{N}_{r_k}$. That is $r_{10 \times k + 1} = r_{10 \times k + 2} = \ldots = r_{10 \times k + 10}$ for $\forall k \in \mathbb{N}$. $\mathcal{Q}_1$ and $\mathcal{Q}_2$  are defined accordingly.
\end{itemize}   

\begin{figure}
    \centering
    \includegraphics[width=\textwidth]{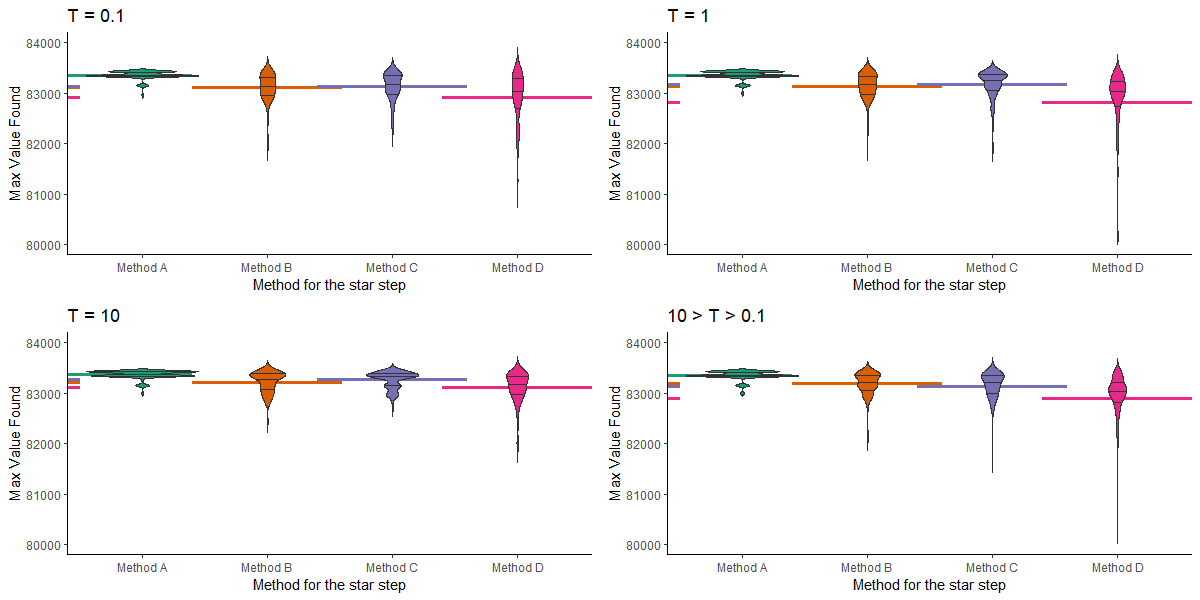} \\
    \caption{Comparison of different methods to choose the subsets for PNS, in terms of the highest (log) target density value $\log \pi(x) = x^T Q x$ found. Method A: random subset every step; method B: random subset every ten steps; method C: systematic subset every step; method D: systematic subset every ten steps. Random upper triangular QUBO matrix where the non-zero elements are generated by $Q_{i, j} \sim N(0, 100^2)$. Four different cooling schedules where T = 0.1, 1, and 10 for all $n$, and T being geometric from 10 to 0.1, are used here. The number of iterations for all methods is 1000. The three black lines inside the violin plots are $25\%$, $50\%$, and $75\%$ quantile lines. The colored segments represent the mean values.}
    \label{fig-howtochoose}
\end{figure}

Again, we use the $200\times200$ QUBO example. The settings for the simulation are the same as in Section~\ref{sec-qubo}. For Method C and D, the two Partial Neighbor Sets $ \mathcal{N}_1$ and $\mathcal{N}_2$ are defined to be flipping the first 100 entries in $x$ and flipping the last 100 entries in $x$. The result for the simulation is shown in Figure~\ref{fig-howtochoose}. This figure shows that the random subset at every step (Method A) performs the best in all four Cooling Schedules. Therefore, we will keep using Method A in all later parts.

\begin{figure}
    \centering
    \includegraphics[width=\textwidth]{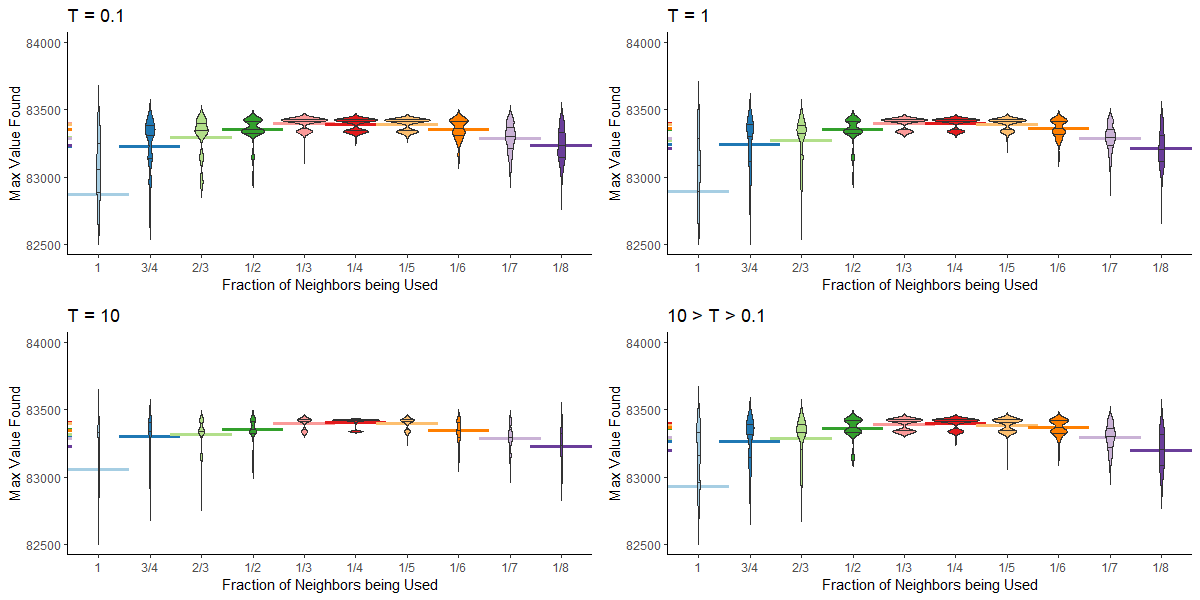} \\
    \caption{Comparison of different sizes of the random subsets for PNS, in terms of the highest (log) target density value $\log \pi(x) = x^T Q x$ being found. Subset sizes are $N \times \{1, \frac{3}{4}, \frac{2}{3}, \frac{1}{2}, \frac{1}{3}, \frac{1}{4}, \frac{1}{5}, \frac{1}{6}, \frac{1}{7}, \frac{1}{8}\}$. Random upper triangular QUBO matrix where the non-zero elements are generated by $Q_{i, j} \sim N(0, 100^2)$. Four different cooling schedules where T = 0.1, 1, and 10 for all $n$, and T being geometric from 10 to 0.1, are used here. The number of iterations for all methods is 1000. The three black lines inside the violin plots are $25\%$, $50\%$, and $75\%$ quantile lines. The colored segments represent the mean values.}
    \label{fig-howmuchchoose}
\end{figure}

In addition, we used Partial Neighbor Sets with half elements from all neighbors in previous simulations. Now we compare the Partial Neighbor Sets with cardinality of $\lvert \mathcal{N}(X) \rvert \times \{1, \frac{3}{4}, \frac{2}{3}, \frac{1}{2}, \frac{1}{3}, \frac{1}{4}, \frac{1}{5}, \frac{1}{6}, \frac{1}{7}, \frac{1}{8}\}$ by the same simulation settings as before. From Figure~\ref{fig-howmuchchoose}, we can see that $\frac{1}{3}, \frac{1}{4}, \frac{1}{5}$ are overall the best among all the choices. Thus, we can conclude that Partial Neighbor Sets with around $25\%$ of the neighbors being considered at each step are the best for the QUBO question stated above.

Therefore, we conclude that our best method to do optimization for the $200\times200$ QUBO question is Algorithm~\ref{alg-qubo}.

\begin{algorithm}
\caption{Partial Neighbor Search for the 200 by 200 QUBO question}\label{alg-qubo}
\begin{algorithmic}
\State initialize $J_0$
\For{$k$ in $1$ to $K$}
    \State randomly pick $ \mathcal{N}_k(J_{k-1}) \subset  \mathcal{N}(J_{k-1})$ where $\lvert  \mathcal{N}_k(J_{k-1}) \rvert = 50$ 
    \State \Comment{Only 50 out of the 200 neighbors will be considered}
    \For{$Y \in \mathcal{N}_k(J_{k-1})$}
        \State calculate $q(Y) = \min\{1, [\frac{\exp(Y^T Q Y)}{\exp(J_{k-1}^T Q J_{k-1})}]^{\frac{1}{T(k)}}\}$
        \State \Comment{the transition prob. from $J_{k-1}$ to $Y$}
    \EndFor
    \State choose $J_{k} \in  \mathcal{N}_k(J_{k-1})$ such that $\hat{P}(J_{k} = Y \mid J_{k-1}) \propto q(Y)$
    \State \Comment{choose the next Jump Chain State}
\EndFor
\end{algorithmic}
\end{algorithm}

\section{Comparison with Tabu Rejection-Free algorithm} \label{sec-tabu}

Tabu search \citep{glover1989tabu} \citep{glover1990tabu} is also a methodology in optimization that guides a local heuristic search procedure to explore the solution space beyond local optimality. The idea of Tabu search is to prohibit access to specific previously-visited solutions. Tabu search is the most intuitive method to help the Markov Chain escape from local maximum areas, as in Figure~\ref{fig-stuckexample}. After moving from state A to state B, we must choose our next state among $\{B_1, B_2, \dots, B_N\}$. We can combine our Rejection-Free algorithm for optimization with Tabu search and then compare this new method to the PNS by the QUBO question. Note that we do not need to record all visited states since we are almost impossible to revisit a state after a certain number of steps. Thus, we only need to record the last several steps and prohibit our Markov chain from revisiting them. The new algorithm is formulated as Algorithm \ref{alg-tabu}.

\begin{algorithm}
\caption{L steps Simplified Tabu Rejection-Free for optimization}\label{alg-tabu}
\begin{algorithmic}
\State initialize $J_0$
\For{$k$ in $1$ to $K$}
    \For{$Y \in \mathcal{N}(J_{k-1}) \backslash \{J_{k-2}, \dots, J_{k-L-1}\}$} 
        \State \Comment{Remove states from the last L steps}
        \State $q(Y) = \min\{1, [\frac{\exp(Y^T Q Y)}{\exp(J_k^T Q J_k)}]^{\frac{1}{T(k)}}\}$ 
        \State \Comment{the transition prob. from $J_{k-1}$ to $Y$}
    \EndFor
    \State choose $J_{k} \in  \mathcal{N}_k(J_{k-1})$ such that $\hat{P}(J_{k} = Y \mid J_{k-1}) \propto q(Y)$
    \State \Comment{choose the next Jump Chain State}
\EndFor
\end{algorithmic}
\end{algorithm}

\begin{figure}
    \centering
    \includegraphics[width=\textwidth]{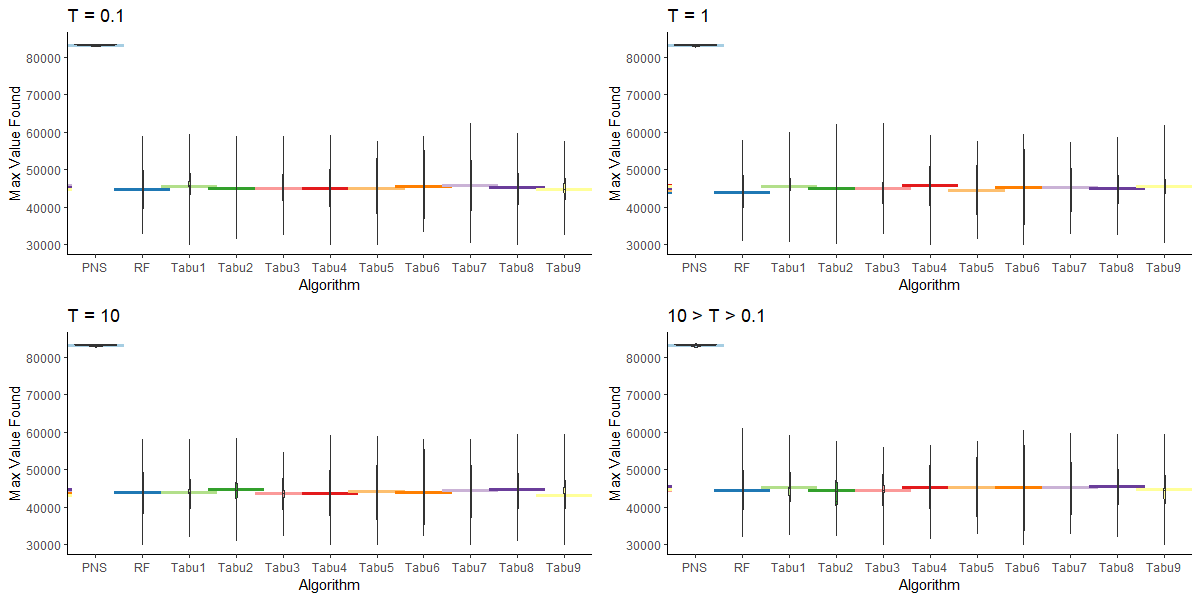}
    \caption{Comparison of PNS, Rejection-Free and 1-Step to 9-steps Simplified Tabu Rejection-Free, in terms of the highest (log) target density value $\log \pi(x) = x^T Q x$ found. Random upper triangular QUBO matrix where the non-zero elements are generated by $Q_{i, j} \sim N(0, 100^2)$. Four different cooling schedules where $T = 0.1$, $1$, and $10$ constantly, and T being geometric from 10 to 0.1, are used here. The run time for all algorithms on a single-core implementation is about the same. The number of iterations for PNS is 400, and the number of iterations for all other methods is 100. The colored segments represent the mean values.}
    \label{fig_compare_tabu}
\end{figure}

Here, we compare PNS with L-step Simplified Tabu Rejection-Free for $L = 1, 2, 3, \dots, 9$. Again, we randomly generate a 200 by 200 upper triangular QUBO matrix. The non-zero elements from the 200 by 200 upper triangular matrix $Q$ were generated randomly with $Q_{i, j} \sim N(0, 100^2)$ for $i < j$. Note that we need to consider about $200$ neighbors at each step for both Rejection-Free and Simplified Tabu Rejection-Free, while we only need to consider $50$ neighbors at each iteration for PNS. If we proceed with the algorithms with a single-core implementation, Rejection-Free and Tabu Rejection-Free need about four times longer than PNS with the same number of steps. Therefore, we can compare the PNS with $4 \times 100 = 400$ iterations with the other methods to get a fair comparison for the program on a single core. Note that we are using this many numbers of steps here because 400 steps are enough for PNS to find a good enough answer. The result for the simulation is shown in Figure~\ref{fig_compare_tabu}. From this plot, we can see that PNS performs much better than Rejection-Free and Simplified Tabu Rejection-Free.

\section{Application to Knapsack problem} \label{sec-knapsack}

The Knapsack problem is another well-known NP-hard problem in optimization \citep{salkin1975knapsack}. We consider the simplest 0-1 Knapsack problem here. Given a knapsack of max capacity $W$ and $N$ items with corresponding values $\{v_i\}_{i=1}^N$ and weights $\{w_i\}_{i=1}^N$, we want to find a finite number of items among all $N$ items which can maximize the total value while not exceeding the max capacity of the knapsack. That is, for given $W > 0$, $\{v_i\}_{i=1}^N > 0$ and $\{w_i\}_{i=1}^N > 0$, find a sequence of $N$ binary variable $\{X_i\}_{i=1}^N \in \{0, 1\}$ to maximize

    \begin{equation}
    \begin{aligned}
        & \sum_{i=1}^N v_i X_i \\
        \text{subject to } & \sum_{i=1}^N w_i X_i \le W 
    \end{aligned}
    \end{equation}
 
Since the Knapsack problem is NP-hard, we can use the Simulated Annealing algorithm to find a feasible solution. For this simulation, we set $W = 100,000$. We randomly generate $N = 1000$ items where the values and weights are random by $w_i, v_i \sim \mbox{Poisson}(1000)$. The mean and the variance for $\mbox{Poisson}(1000)$ are both 1000. Suppose we want to find a binary vector $X = (X_1, X_2, \dots, X_N)^T$ of dimension $N$ to maximize $v^T X$ subject to $w^T X \le W$. 

Again, we used a uniform proposal distribution among all neighbors where the neighbors are defined as binary vectors with Hamming distance 1. That is, $\mathcal{Q}(X, Y) = \frac{1}{N}$ for $\forall Y \in \mathcal{N}(X)$, where $Y \in \mathcal{N}(X) \iff \lvert X - Y \rvert = \sum_{i=1}^N \lvert X_i - Y_i \rvert = 1$, $\forall X, Y \in \{0, 1\}^N$. We randomly choose half of the neighbors at each step for PNS. That is, $\lvert \mathcal{N}_k(X) \rvert = \frac{1}{2} \lvert \mathcal{N}(X) \rvert = 500$ for $\forall X \in \{0, 1\}^N$. Moreover, the target density $\pi(X) = \mathbb{1}(w^T X \le W) \times v^T X$, where $\mathbb{1}$ represents the indicator function. In addition, $T(k)$ represents the temperature at step $k$ for the Cooling Schedule here.

\begin{figure}
    \centering
    \includegraphics[width=\textwidth]{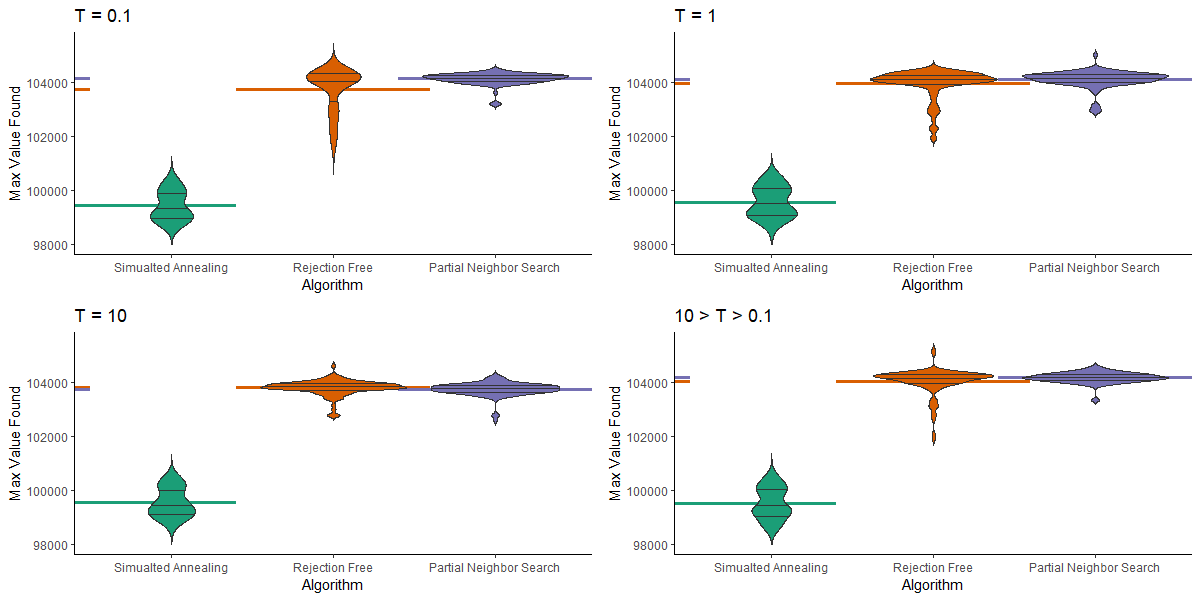}
    \caption{Comparison of Simulated Annealing, Rejection-Free, and PNS in terms of the highest target density values found in Knapsack Problem with $W = 100,000$, $N = 1000$, $w_i, v_i \sim \mbox{Poisson}(1000)$. Four different cooling schedules where $T = 0.1$, $1$, and $10$ constantly, and T being geometric from 10 to 0.1, are used there. The number of iterations for Simulated Annealing is 1,000,000, while the number for Rejection-Free and PNS is 1000. The three black lines inside the violin plots are $25\%$, $50\%$, and $75\%$ quantile lines. The colored segments represent the mean values.}
    \label{fig_knspsack}
\end{figure}

Again, we compare Simulated Annealing, Rejection-Free with PNS here. The result is shown in Figure~\ref{fig_knspsack}. The plot shows that Rejection-Free for optimization and PNS algorithm are better than the regular Simulated Annealing algorithm in all four Cooling Schedules. Again, for the simulation shown in Figure~\ref{fig_knspsack}, the numbers of iterations used for the three methods are set to be different to have a fair comparison between three methods. We set the number of iterations for Simulated Annealing to be $1000,000$. The numbers of iterations for Rejection-Free and PNS are $1000$ since we need to consider $1000$ neighbors at each iteration for Rejection-Free for optimization. In contrast, we only need to consider one neighbor for each iteration in Simulated Annealing. 

This result shows that PNS is not always that much better than Rejection-Free when the number of iterations is the same. In some cases, where the target distribution is not sharply peaked, and there are not too many local extreme areas, Rejection-Free can also have excellent performance. Note that if we run the above simulation on a single core, PNS will only take about half of the time used by Rejection-Free, and if we use parallel hardware to apply the above algorithm, Rejection-Free and PNS will take about the same time.

In addition, Rejection-Free is not always better than simple Simulated Annealing. For example, if $\pi(X) \equiv 1$ for all $X \in S$, there will be no rejections. The Simulated Annealing will move to a new state by computing a single probability, while the Rejection-Free will do the same but compute the probabilities for all neighbors. However, when the dimension of the problem is large, or the target density is sharply peaked, the PNS will perform much better than Rejection-Free, and Rejection-Free will perform much better than Simulated Annealing.

\section{Application to 3R3XOR problem} \label{sec-xor}

The 3R3XOR problem is a methodology for generating benchmark problem sets for Ising machines devices designed to solve discrete optimization problems cast as Ising models introduced by \cite{hen2019equation}. The Ising model, named after Ernst Ising, is concerned with the physics of magnetic-driven phase transitions \citep{cipra1987introduction}. The Ising model is defined on a lattice, where a spin $s_i \in \{-1, 1\}$ is located on each lattice site \citep{block2012computer}. The optimization question for the Ising model has been widely applied to many scientific problems such as neuroscience \citep{hopfield1982neural} and environmental science \citep{ma2014ising}. Thus, algorithms, even special-purpose programmable devices, designed to solve discrete optimization problems cast as Ising models are popular \citep{hen2019equation}, and our PNS algorithm is one of them. 

However, the non-planar Ising model is NP-complete \citep{cipra2000ising}. We cannot find an optimal state from an Ising model in polynomial time. Then, it is hard for us to compare the performance of the heuristic solvers, such as Rejection and PNS, by the time used to find the optimal state from a random Ising model. On the other hand, \cite{hen2019equation} introduced a tool for benchmarking Ising machines in 2019. In his approach, linear systems of equations are cast as Ising cost functions. The linear systems can be solved quickly, while the corresponding Ising model exhibits the features of NP-hardness \citep{hen2019equation}. This way, we can construct special Ising models with a unique known optimal state. Then we can use these special Ising models to compare the heuristic solvers' runtimes for finding the optimal state.

In this section, we focus on constructing a simplified version of 3-body Ising with $N$ spins from a binary linear system of $N$ equations. The simplified version is defined as follows:
\begin{equation}
    H(\{s_j\}) = \sum_{a < b < c} \mathbf{M}_{a, b, c} s_a s_b s_c,
\end{equation}
where $s_i\in \{-1, 1\}$ for $\forall i = 1,2,\ldots, N$. $\mathbf{M}_{a, b, c}$ is a $N \times N \times N$ matrix where $\mathbf{M}_{a, b, c} = 0$ $\forall a \ge b$, $b \ge c$, or $a \ge c$. 

In \citeauthor{hen2019equation}'s (\citeyear{hen2019equation}) approach, we start by choosing a binary matrix $\{ \mathbf{A}_{i, j}\}$ and a binary vector $\{ \mathbf{b}_j\}$ to form a modulo 2 linear system of N equations in N variables.
\begin{equation}
    \sum_{j=1}^N  \mathbf{A}_{i, j} x_j \equiv \mathbf{b}_i \mod 2 \mbox{, for } i = 1,2,\dots,N.
\end{equation}
This module 2 linear system of equations can always be solved in polynomial time using Gaussian elimination. In addition, as long as the binary matrix $\{ \mathbf{A}_{i, j}\}$ is invertible, the solution (if exists) is unique. Suppose $\{x_1, . . . , x_n\}$ are $n$ binary variables. Then for given $\{ \mathbf{A}_{i, j}\}$ and $\{ \mathbf{b}_j\}$, we define 
\begin{equation}
    F(\{x_j\}) = \sum_{i=1}^N \mathbb{1}\Big(\sum_{j=1}^N   \mathbf{A}_{i, j} x_j \not\equiv  \mathbf{b}_i \mod 2\Big),
\end{equation}
where $\mathbb{1}$ means indicator function here. Since $F$ is a sum of $N$ indicator functions, then $0 \le F \le N$ and the minimum bound is reached when $\{x_j\}$ is the solution to the modulo 2 linear system.

Let $s_j = 1 - 2 x_j \in \{-1, 1\}\mbox{ for } j = 1, 2, \dots, N$ be $N$ Ising spins. Then we must have
\begin{equation}
    \prod_{j: \mathbf{A}_{i, j} = 1} s_j = (-1)^{ \mathbf{b}_i} \mbox{ if and only if } \sum_{j=1}^N   \mathbf{A}_{i, j} x_j \equiv  \mathbf{b}_i \mod 2,
\end{equation}
$\forall i = 1,2,\dots,m$. Then 
\begin{equation}
\begin{aligned}
    F & = \sum_{i=1}^N \mathbb{1}\Big(\sum_{j=1}^N   \mathbf{A}_{i, j} x_j \not\equiv  \mathbf{b}_i \mod 2\Big) \\
    & = \sum_{i=1}^N \mathbb{1}\Big( \prod_{j:  \mathbf{A}_{i, j} = 1} s_j \ne (-1)^{ \mathbf{b}_i} \Big)
    \mbox{, since} \prod_{j:  \mathbf{A}_{i, j} = 1} s_j \mbox{ and } (-1)^{ \mathbf{b}_i} \in \{-1, 1\}\\
    & = \frac{1}{2}\Big[\sum_{i=1}^N \Big(1 -  (-1)^{ \mathbf{b}_i} \prod_{j:  \mathbf{A}_{i, j} = 1} s_j \Big)\Big].
\end{aligned}
\end{equation}
After dropping immaterial constants, we define 
\begin{equation}
    F_0(\{s_j\}) = \sum_{i=1}^N \Bigg[(-1)^{ \mathbf{b}_i} \prod_{j:  \mathbf{A}_{i, j} = 1}  s_j \Bigg].
\end{equation}
Note that $F \ge 0$ and the minimum bound is reached when $\{x_j\}$ is the solution to the modulo 2 linear system. Thus, $F_0 \le N$, and the maximum bound will be reached when $\{x_j \mid x_j = \frac{1}{2}(1 - s_j)\}$ is the solution to the modulo 2 linear system. In addition, as long as the matrix $\{\mathbf{A}_{i, j}\}$ is invertible, the solution to the equation system must uniquely exist, and then there must exist a single configuration maximize $F_0$ whose maximum value is exactly $N$. 

Again, the Hamiltonian for simplified 3-body Ising model including only the cubic term to be $H(\{s_j\}) = \sum_{a < b < c} \mathbf{M}_{a, b, c} s_a s_b s_c$. Here, we assume, on each row of binary matrix $\{\mathbf{A}_{i, j}\}$, $\sum_{j=1}^N \mathbf{A}_{i, j} = 3$. Then, let $\mathbf{M}_{a, b, c} = (-1)^{\mathbf{b}_i}$ if $\exists i, a < b < c$ such that $\mathbf{A}_{i, a} = \mathbf{A}_{i, b} = \mathbf{A}_{i, c} = 1$, and $\mathbf{M}_{a, b, c} = 0$ otherwise. Then, we have $H(\{s_j\}) = F_0(\{s_j\})$.

Thus, we can construct an Ising model with a unique optimal bound with a known optimal value $N$ as follows:

\begin{enumerate}
    \item find an invertible binary matrix $\{\mathbf{A}_{i, j}\}$ and a binary vector $\{\mathbf{b}_i\}$, where $\sum_{j=1}^N \mathbf{A}_{i, j} = 3$, $\forall i$
    \item solve the modulo 2 linear equation system $\sum_{j=1}^N \mathbf{A}_{i, j} x_j \equiv \mathbf{b}_i \mod 2$, for $i = 1,2,\dots,N $ to make sure the unique solution exists
    \item define $\mathbf{M}_{a, b, c}$ be a $N \times N \times N$ matrix where  $\mathbf{M}_{a, b, c} = (-1)^{\mathbf{b}_i}$ if $\exists i, a < b < c$  such that $\mathbf{A}_{i, a} = \mathbf{A}_{i, b} = \mathbf{A}_{i, c} = 1$, and $\mathbf{M}_{a, b, c} = 0$ otherwise
    \item then we must have a unique optimal solution $s_{\max}$ for $H(s_{\max}) = \max(H(s)) = N$
\end{enumerate}

\begin{figure}
    \centering
    \includegraphics[width= \linewidth]{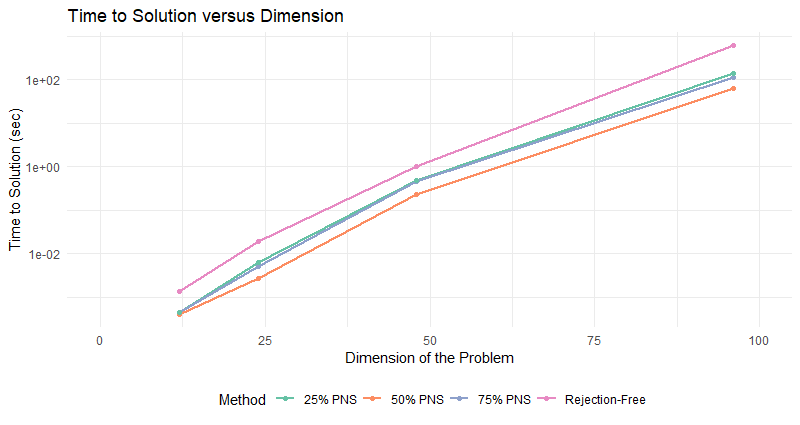}
    \caption{Comparison of the minimum value for the time used to find the optimal state by Rejection-Free and PNS with $25\%$, $50\%$, and $75\%$ of the neighbors being considered at each step for a random Ising model generated by 3R3XOR. Each dot represents the median of 50 repeated simulations for a given problem size $N = 12$, $24$, $48$ and $96$.}
    \label{fig_3r3xor}
\end{figure}

By constructing the special 3-body $N \times N \times N$ Ising model with a unique optimal solution of maximum bound $N$, we can examine the performance of the Rejection-Free and PNS algorithms on these special Ising models. Again, uniform proposal distributions are used here, and the neighbors are defined as binary vectors with Hamming distance 1. We random generate the special Ising models with four different sizes $N = 12$, $24$, $48$ and $96$. For each of these four different sizes, we generate 50 different Ising models and record the time used by the algorithms to reach the unique optimal state. The median of these 50 results for both Rejection-Free and PNS algorithms are shown in Figure~\ref{fig_3r3xor}. From this figure, Rejection-Free is the worst. $25\%$ PNS performs comparably to $75\%$, and the $50\%$ PNS performs the best.

\section{Application to Continuous State Space} \label{sec-continuous}

In previous sections, we focused on optimization questions with the discrete state space $\mathcal{S}$ where all states have at most a finite number of neighbors. Meanwhile, Simulated Annealing works for general state space. In addition, Theorem $13$ in \cite{rosenthal2021jump} extended the Rejection-Free for sampling to general state space. Similarly, we can extend the Rejection-Free for optimization to general state space. 

Although we have a solid theory base for Rejection-Free in general state space, it is challenging to apply Rejection-Free to those cases. There is a major difficulty involved in the for loop that calculates the transition probability of all neighbors in Algorithm~\ref{alg-rf}. In continuous cases, although numerical integration of all transition probability can be performed, it is unlikely that such tasks may be efficiently divided among specialized hardware with a certain number of parallel processing units. On the other hand, PNS, as described in Algorithm~\ref{alg-pns}, can be applied straightforwardly to continuous cases by choosing the Partial Neighbors Sets $\mathcal{N}_k(X)$ to be finite subsets of all the neighbors $\mathcal{N}(X)$ in Algorithm~\ref{alg-pns}.

We compare the performance of Simulated Annealing with our PNS on a simple example of quadratic programming, which belongs to the category of continuous optimization, as stated below:
\begin{equation}
\begin{aligned}
    \arg \max \mbox{ } & x^T Q x \\
    \mbox{subject to } & x_i \ge 0 \mbox{, } \forall i = 1, 2, \dots, N \\
    & \sum_{i=1}^N x_i = 1,
\end{aligned}
\end{equation}
where $Q$ is a given an upper triangular $N$ by $N$ matrix and $x \in \mathbb{R}^N$. For most cases, the quadratic programming is stated by $\arg\min$ instead of $\arg\max$. We use the $\arg\max$ version here to be consistent with the QUBO question in Section~\ref{sec-qubo}, and $\arg\max$ is equivalent to $\arg\min$ when replacing $Q$ by $-Q$. This quadratic programming question is also NP-hard as long as $Q$ is indefinite \citep{sahni1974computationally}, where indefinite means matrices that are neither positive semi-definite nor negative semi-definite. 

We randomly generate a 200 by 200 upper triangular to be the matrix $Q$, where the non-zero elements from the 200 by 200 upper triangular matrix $Q$ were generated randomly by $Q_{i, j} \sim \text{Normal}(0, 100^2)\mbox{, } \forall i \le j$. We compare Simulated Annealing and PNS in 100 simulation runs here. We omit Rejection-Free since applying Rejection-Free to continuous cases is quite hard. 

The target density value is set to be $\pi(x) = \exp\{x^T Q x\}$, $\forall x$ such that $x_i \in (0, 1)$, $\forall i = 1,2,\dots,N$, and $\pi(x) = -\infty$ otherwise. In addition, the proposal distribution $\mathcal{Q}$ and the corresponding neighbor set $\mathcal{N}$ are defined as follows:
\begin{enumerate}
    \item for state $x = (x_1, x_2, \dots, x_N)^T \in \mathcal{S}$, choose a random entry $x_r$ for $r \in \{1, 2,\dots, N\}$;
    \item generate a random value $s \sim \text{Normal}(0, 0.1^2)$;
    \item let $y_r = x_r + s$ and $y_n = x_n \times \frac{1 - x_r}{1 - y_r}$, $\forall n \ne r$;
    \item if $y_r \notin (0, 1)$, then the corresponding $\pi(y)$ is defined to be $-\infty$; in practice, we just need to generate a new $y$; also note that, as long as $y_r, x \in (0, 1)$, we must have $y \in (0, 1)$ as well;
    \item to ensure the reversibility within each Partial Neighbor Set, we also consider $y'_r = x_r - s$ and $y'_n = x_n \times \frac{1 - x_r}{1 - y'_r}$, $\forall n \ne r$; if $y'_r \notin (0, 1)$, then we can ignore $y'$.
\end{enumerate}
With the given steps, we have $\sum_{n=1}^N y_n = 1$ as long as $\sum_{n=1}^N x_n = 1$. This method is similar to component-wise Simulated Annealing. We find a random component, magnify or minify it, and then modify the rest of the entries accordingly to make the summation remain unchanged. This proposal distribution $\mathcal{Q}$ is therefore systematic. By the above ways to generate neighbors, we can eliminate the constraints that $x_i \ge 0 \mbox{, } \forall i = 1, 2, \dots, N$, and $\sum_{i=1}^N x_i = 1$, and we only need to focus on $\arg \max x^T Q x$. 

For Simulated Annealing, we randomly generate one neighbor by the above given steps and calculate the transition probability. For PNS, we can generate, for example, $20$ random neighbors at each step. In this case, the Partial Neighbor Set $\mathcal{N}_i$ is only a random subset of $\mathcal{N}$ with $20$ elements, and thus, the implementation of PNS is simple compared to the Rejection-Free. 

\begin{figure}
    \centering
    \includegraphics[width=\textwidth]{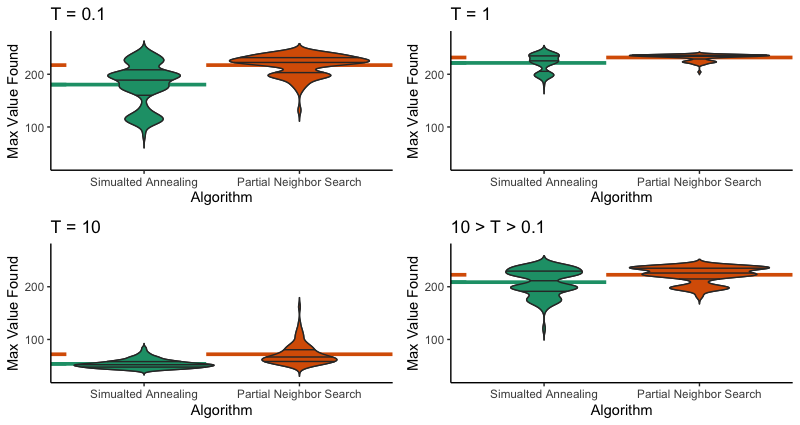} \\
    \caption{Comparison of Simulated Annealing and PNS in terms of the highest (log) target distribution value $\log \pi(x) = x^T Q x$ being found, for a random upper triangular matrix $Q$ and $x \in \mathbb{R}^N$ subject to $x_i \ge 0$, $\forall i = 1, 2, \dots, N$, and $\sum_{i=1}^N x_i = 1$. The non-zero elements are generated by $Q_{i, j} \sim N(0, 100^2)$. Four different cooling schedules where $T(k) = 0.1$, $1$ and $10$ constantly, and $T(k)$ being geometric from 10 to 0.1 are used here. The number of iterations for Simulated Annealing is $600,000$, and the number of iterations for PNS is $72,000$. The run times for these two algorithms on a single-core implementation are both around 80 seconds. The three black lines inside the violin plots are $25\%$, $50\%$, and $75\%$ quantile lines. The colored segments represent the mean values.}
    \label{fig_continuous}
\end{figure}

The result for the simulation is shown in Figure~\ref{fig_continuous}. We can see that the PNS performs better than Simulated Annealing in all four different cooling schedules. However, the difference between PNS and Simulated Annealing in this continuous example is not as much as the difference between the algorithms from the discrete QUBO questions. This is because the continuous example is not as sharply peaked as the discrete example from Section \ref{alg-qubo}. After we choose a random entry $r$, we only need to move a small step around the original value of $x_r$. On the other hand, we have to flip between $0$ and $1$ in the discrete example. Thus, the rejection rate for the Simulated Annealing is lower than the rate from the discrete example, so the performance of these two algorithms gets closer.

In addition, PNS is specially designed for parallelism hardware. Again, with a specialized dedicated processor such as DAU, PNS can yield 100x to 10,000x speedups \cite{sheikholeslami2021power}. In addition, this example also shows PNS has more flexibility compared to the Rejection-Free algorithm. Again, Rejection-Free can hardly work for cases with infinite neighbors, while PNS can be easily applied by choosing finite $\mathcal{N}_k$. 

Moreover, the number of elements in $\mathcal{N}_k$ needs to be reasonable for PNS to keep its performance. For example, if we used $\lvert \mathcal{N}_k \rvert = 500$, we would calculate too many transition probabilities at each step, and the algorithm would be inefficient. Meanwhile, if we used  $\lvert \mathcal{N}_k \rvert = 2$, the number of Partial Neighbor Sets being considered at each step would be too few. As PNS will force the Markov chain to move to one element from the Partial Neighbor Set $\mathcal{N}_k$, it will move to some terrible choices of states when all states in the Partial Neighbor Set $\mathcal{N}_k$ have small target distribution values. In the above simulation, choosing $\lvert \mathcal{N}_k \rvert$ from $10$ to $30$ won't make a big difference.

\section{Summary} \label{sec-summary}

This paper illustrates Rejection-Free Simulated Annealing algorithms that consider all neighbors at each step in order to prevent inefficiency from rejections. We have also proposed a Partial Neighbor Search (PNS) algorithm based on the Rejection-Free technique in order to address the issue of local maximum area. Three sets of discrete examples have been simulated to demonstrate that PNS can produce significant speedups in optimization problems. PNS has also been applied to continuous examples in order to demonstrate its greater flexibility in comparison to Rejection-Free. 

\section*{Acknowledgments}

The author(s) would like to thank Fujitsu Ltd. and Fujitsu Consulting (Canada) Inc. for providing financial support.

\bigskip

\bibliography{References}

\end{document}